\documentclass{article}

\usepackage{iclr2027_conference,times}
\iclrfinalcopy

\usepackage{amsmath,amssymb,graphicx}
\usepackage{hyperref}
\usepackage{url}
\usepackage{algorithm}
\usepackage{algorithmic}
\usepackage{booktabs}
\usepackage{caption}
\usepackage{xurl}
\usepackage{pdfpages}
\usepackage{float}
\usepackage{tabularx}

\title{Damped Gauss–Newton Search for Multi-Metric Hyperparameter Optimization}

\author{Qinwu Xu, Zhanyu Ge, Yifan Jiang}

\begin{document}
\maketitle
\fancyhead{} %to remove heads

\begin{abstract}
We study hyperparameter optimization (HPO) from a numerical-optimization perspective and propose a black-box, target-seeking damped Gauss--Newton method for multiple validation responses. In contrast to uninformed candidate search such as grid or random search, the method uses numerical sensitivity information inferred from previously observed changes in hyperparameters and validation performance. Successive full-vector observations are used to construct an iterative secant approximation of the local Jacobian, which relates changes in the performance vector to changes in all optimized hyperparameters; a damped Gauss--Newton system then produces a joint update of the complete hyperparameter vector rather than optimizing coordinates independently. The method does not differentiate through model training or require dedicated coordinate-wise perturbation runs: following the initial two evaluations used to establish the first secant,
each subsequent iteration requires one new full-vector evaluation. Tikhonov regularization stabilizes the resulting local inverse problem when the number of optimized parameters exceeds the number of performance responses. Unlike conventional multi-objective HPO aimed at approximating a Pareto front, our formulation seeks a specified target operating point in the multi-response performance space. We evaluate the same optimization principle in three complementary settings: four-dimensional XGBoost training-time HPO on three public classification datasets, six-dimensional LoRA/SFT post-training of Qwen3-VL-8B on a binary VQAv2 subset, and eight-dimensional post-processing threshold optimization on fixed classifier outputs. On the XGBoost benchmarks, the proposed method achieves validation performance comparable to grid search, random search, and TPE under a smaller number of model evaluations, although it does not consistently outperform these black-box baselines. The VLM study demonstrates joint post-training hyperparameter updates using TPR/FPR feedback and reveals strongly non-monotonic response and sensitivity trajectories, including a high-TPR/high-FPR degenerate regime. The threshold study further demonstrates target seeking in an underdetermined two-response, eight-parameter system and shows sensitivity to initialization and the importance of retaining the best observed iterate. Across these settings, the results support trajectory-based secant sensitivity as an evaluation-efficient local optimization mechanism complementary to global black-box HPO.
\end{abstract}

\section{Introduction}
\label{sec:introduction}

Hyperparameter optimization (HPO) is a central component of modern
machine learning because choices such as learning rate, regularization,
subsampling, model depth, and optimization schedules can substantially
affect predictive performance, training stability, and computational
cost~\cite{bischl2023hyperparameter}. Unlike model parameters learned
directly through training, hyperparameters are typically selected by an
outer optimization process in which candidate configurations are trained
or evaluated and their validation performance is observed. When each
candidate requires a complete training or post-training run, the number
of such evaluations can dominate the cost of HPO.

Classical search methods provide simple but largely direction-uninformed
strategies. Grid search evaluates a predefined Cartesian product of
candidate values, with cost growing rapidly with the number of
hyperparameters. Random search can be substantially more efficient by
exploring more distinct values along influential dimensions
~\cite{bergstra2012random}, but successive candidates do not explicitly
exploit an estimated local response direction derived from previous
performance changes.

Sequential model-based methods use previous evaluations more explicitly.
Bayesian optimization constructs a probabilistic surrogate of validation
performance and uses an acquisition function to select subsequent
configurations~\cite{snoek2012practical}. TPE uses density estimation for
conditional and mixed search spaces~\cite{bergstra2013making}, while
SMAC provides a widely used sequential model-based framework for
algorithm configuration and HPO
~\cite{hutter2011sequential,lindauer2022smac3}. These approaches exploit
information from previous trials through a learned surrogate or density
model. The approach studied here instead extracts a local numerical
sensitivity directly from changes observed along the optimization
trajectory and uses it to construct the next joint hyperparameter update.

Other HPO methods reduce cost through adaptive resource allocation or
population-based exploration. Hyperband terminates poorly performing
configurations early~\cite{li2017hyperband}, BOHB combines model-based
search with Hyperband~\cite{falkner2018bohb}, and population-based
training jointly trains and adapts a population of models
~\cite{jaderberg2017population}. Evolutionary strategies provide another
derivative-free alternative, including covariance-adaptation methods that
update a search distribution from populations of evaluated solutions
~\cite{hansen2001cma}. These approaches are particularly useful when
partial training trajectories or parallel evaluations are available,
but typically explore multiple candidates or maintain a population or
distribution over candidate configurations.

Gradient-based HPO provides a different source of directional
information. By expressing training and hyperparameter selection as a
bilevel optimization problem, hypergradients can be obtained by
differentiating through, approximating, or implicitly differentiating
the inner training procedure~\cite{franceschi2018bilevel}. Reversible
differentiation, for example, propagates gradients through the training
trajectory and can optimize many hyperparameters simultaneously
~\cite{maclaurin2015gradient}. Such methods provide rich gradient
information when the learning pipeline is differentiable and its
optimization dynamics are accessible. These requirements can limit
direct application to black-box learners, non-differentiable
hyperparameters, legacy training systems, or evaluation metrics for
which differentiation through the complete training procedure is
unavailable.

Our approach instead obtains directional information numerically while
keeping the learning procedure itself as a black box. It is therefore
related to derivative-free numerical optimization, where local models,
finite differences, or function-value observations are used when
analytical derivatives are unavailable
~\cite{conn2009derivative,khare2021trust,fortin2004trust}. SPSA, for
example, estimates a search direction from simultaneous random
perturbations rather than coordinate-wise finite differences
~\cite{spall1992multivariate}. In contrast, our method introduces no
dedicated perturbation evaluations to estimate a scalar-objective
gradient. Successive \emph{full-vector} HPO iterates instead provide
paired changes in the hyperparameter and validation-performance vectors.
These trajectory observations are reused as local secant sensitivity
information, from which a regularized Gauss--Newton system produces a
joint update of the complete hyperparameter vector. Thus, the method is
sensitivity-informed without requiring analytical gradients or
backpropagation through model training.

Multiple validation responses introduce an additional distinction.
Multi-objective HPO (MOHPO) is commonly formulated as the search for a
set of nondominated configurations representing tradeoffs among competing
objectives, using approaches such as Pareto-based evolutionary algorithms
or surrogate-based methods such as ParEGO
~\cite{deb2002nsga,knowles2006parego,morales2023survey,karl2022multiobjective}.
We consider a related but different \emph{target-seeking} problem. Let
$\boldsymbol{\phi}\in\mathbb{R}^{d}$ denote the hyperparameter vector and
$\mathbf{P}(\boldsymbol{\phi})\in\mathbb{R}^{m}$ the corresponding vector
of measured validation responses. 

Given a desired operating point
$\mathbf{P}^{\star}\in\mathbb{R}^m$, the target-seeking problem is
\begin{equation}
\min_{\boldsymbol{\phi}}
\frac{1}{2}
\left\|
\mathbf{P}^{\star}-\mathbf{P}(\boldsymbol{\phi})
\right\|_2^2.
\label{eq:intro_objective}
\end{equation}
The objective is therefore not to approximate a Pareto front, but to
move the measured response vector toward a specified operating point
while retaining the individual responses explicitly in the local
numerical system.

At iteration $k$, the black-box learning or evaluation procedure is
executed at the complete vector $\boldsymbol{\phi}_k$ to obtain
$\mathbf{P}_k$. Changes between successive hyperparameter and performance
vectors provide secant information used to update a local sensitivity
matrix. Importantly, this sensitivity should not be interpreted as an
analytical or infinitesimal derivative of the training pipeline.
Validation responses such as classification accuracy, TPR, and FPR can
be noisy, plateaued, or non-monotonic with respect to hyperparameters.
The method instead constructs a finite, trajectory-dependent local
response approximation from observed changes and repeatedly updates it
as new evaluations become available. Damping and best-iterate retention
are consequently important practical components of the procedure.

This construction also has an important evaluation-cost distinction.
A coordinate-wise forward-difference approximation of a
$d$-dimensional sensitivity requires multiple dedicated perturbed
evaluations around a reference configuration. Our iterative construction
instead reuses the sequence of full-vector configurations generated by
the optimization itself. Following the initial two evaluations used to establish the first secant,
each subsequent update introduces one new full-vector evaluation and
jointly advances all optimized hyperparameters. We therefore distinguish the number of
expensive model evaluations from the inexpensive numerical work required
to update the local sensitivity and solve the small regularized linear
system. This distinction is particularly relevant when one evaluation
requires model training or large-model post-training.

When the number of optimized hyperparameters exceeds the number of
performance responses, the local inverse problem is underdetermined.
This occurs in all three studies considered here: four XGBoost
hyperparameters are guided by one primary validation response, six VLM
post-training hyperparameters by two responses, and eight decision
thresholds by two responses. We stabilize these local systems using
Tikhonov regularization in a damped Gauss--Newton update and enforce
feasible parameter ranges where applicable. Because the method is local
and trajectory-based, its behavior can depend on initialization,
damping, and the local response geometry; these properties are examined
empirically rather than assuming global smoothness or convergence.

We evaluate the framework in three complementary settings spanning
training, post-training, and post-processing. First, we jointly tune four
XGBoost training hyperparameters on three public classification datasets
and compare against grid search, random search, and TPE using explicit
candidate-evaluation budgets. Second, we optimize six LoRA/SFT
hyperparameters of Qwen3-VL-8B using TPR and FPR feedback, where each new
HPO evaluation requires VLM post-training and is therefore substantially
more expensive than the tabular case. Third, we optimize eight
class-specific decision thresholds against two target responses using
fixed classifier outputs, isolating the numerical optimization behavior
from model-training cost. Together, these settings test the same
full-vector update principle under substantially different evaluation
loops and response geometries.

Our main contribution is a black-box, target-seeking HPO framework that
reuses successive full-vector evaluations to obtain local secant
sensitivity information and jointly update multiple hyperparameters
through a Tikhonov-regularized damped Gauss--Newton system. The framework
occupies a complementary position between direction-uninformed search,
surrogate-based black-box HPO, and differentiable hypergradient methods:
it extracts numerical directional information from observed HPO
trajectories without differentiating through training, fitting a global
surrogate, or requiring dedicated coordinate-wise perturbation runs.
The experiments characterize both its evaluation efficiency and its
limitations, including non-monotonic responses and sensitivity to
initialization.

% ============================================================
% ============================================================
\section{Method}
\label{sec:method}

\subsection{Problem Formulation}

Let
$\boldsymbol{\phi}=[\phi_1,\ldots,\phi_d]^\top\in\mathbb{R}^d$
denote the vector of parameters to be optimized. For a given
$\boldsymbol{\phi}$, the learning or evaluation procedure produces a
vector of measured validation responses $\mathbf{P}(\boldsymbol{\phi})
=[P_1(\boldsymbol{\phi}),\ldots,P_m(\boldsymbol{\phi})]^\top
\in\mathbb{R}^m .
\label{eq:performance_vector}$ Here, $P_q$ denotes a measured performance quantity rather than an internal training loss; its definition depends on the application.

We consider the black-box mapping $\boldsymbol{\phi}\longmapsto\mathbf{P}(\boldsymbol{\phi}),\label{eq:blackbox}$
without differentiating through the underlying training or evaluation
procedure. Given a desired operating point (performance metrics)
$\mathbf{P}^{\star}\in\mathbb{R}^m$, the target-seeking problem is minimize the difference between it and the predicted performance metrics (equation 1). The individual responses remain explicit rather than being collapsed into a single application-specific scalar score.

The proposed method obtains directional information numerically from
changes in measured performance. In this sense, the search is
\emph{numerical-sensitivity informed}: unlike grid or random search,
which do not construct a local response direction from previous
evaluations, observed parameter--performance changes are used to inform
the next joint update. No analytical hypergradient or differentiation
through model training is required.

This formulation applies to all three settings considered in this paper.
For XGBoost, $\boldsymbol{\phi}$ contains training hyperparameters; for
VLM post-training, it contains LoRA/SFT hyperparameters; and for the
post-processing experiment, it contains class-specific decision
thresholds.

\subsection{Full-Vector Secant Sensitivity}
\label{sec:secant_jacobian}

At iteration $k$, the complete parameter vector
$\boldsymbol{\phi}_k$ is evaluated to obtain $
\mathbf{P}_k=\mathbf{P}(\boldsymbol{\phi}_k).$
Two successive evaluations provide the full-vector secant pair
\begin{equation}
\Delta\boldsymbol{\phi}_k
=
\boldsymbol{\phi}_k-\boldsymbol{\phi}_{k-1},
\qquad
\Delta\mathbf{P}_k
=
\mathbf{P}_k-\mathbf{P}_{k-1}.
\label{eq:iterative_changes}
\end{equation}

We represent the local response through a numerical sensitivity matrix
\begin{equation}
\mathbf{J}_k\in\mathbb{R}^{m\times d},
\qquad
\Delta\mathbf{P}_k
\approx
\mathbf{J}_k\Delta\boldsymbol{\phi}_k .
\label{eq:local_secant}
\end{equation}
The matrix $\mathbf{J}_k$ plays the role of a local Jacobian relating
changes in all $d$ optimized parameters to changes in the $m$ measured
responses. It is not obtained by analytical differentiation through the
learning algorithm; instead, it is a trajectory-dependent numerical
sensitivity estimate constructed and updated from observed full-vector
parameter--performance changes.

Because a successive finite-difference sensitivity estimate requires two
observations, the numerical sensitivity is initialized using a small
simultaneous perturbation of the complete parameter vector. Starting
from the empirical or default initialization $\boldsymbol{\phi}_0$, we
first evaluate
$\mathbf{P}_0=\mathbf{P}(\boldsymbol{\phi}_0)$ and then form $\boldsymbol{\phi}_1
=\boldsymbol{\phi}_0+\delta\mathbf{1}_d,$ where $\mathbf{1}_d\in\mathbb{R}^d$ is the all-ones vector and
$\delta$ is a small perturbation magnitude (e.g., $0.001$), such that
all $d$ optimized coordinates are perturbed simultaneously. A second
full-vector evaluation gives
$\mathbf{P}_1=\mathbf{P}(\boldsymbol{\phi}_1)$, and the resulting
differences
$\Delta\boldsymbol{\phi}_1=\boldsymbol{\phi}_1-\boldsymbol{\phi}_0$
and
$\Delta\mathbf{P}_1=\mathbf{P}_1-\mathbf{P}_0$
provide the initial numerical sensitivity information. Subsequent
sensitivity estimates are obtained directly from consecutive HPO
iterates and therefore require no additional perturbation evaluations.
This initialization requires only one perturbed full-vector evaluation,
rather than separate perturbation evaluations for each optimized
coordinate.

After initialization, the parameter changes generated by the optimizer
itself provide subsequent secant observations. Thus, the numerical
sensitivity is repeatedly updated from the trajectory already produced
by the HPO procedure rather than by introducing dedicated
coordinate-wise perturbation runs at every iteration.

Importantly, $\mathbf{J}_k$ should be interpreted as a finite local
response approximation rather than an analytical or infinitesimal
derivative of the complete learning pipeline. Validation quantities such
as accuracy, TPR, and FPR can be plateaued, noisy, or non-monotonic with
respect to hyperparameters. The sensitivity estimate therefore describes
the locally observed response along the optimization trajectory and is
updated as new observations become available.

For the VLM implementation, the optimized hyperparameters have
substantially different numerical scales. The corresponding secant
quantities are normalized before constructing the numerical sensitivity
used by the controller, preventing parameters with larger numerical
magnitudes from dominating the update solely because of their units.
\subsection{Damped Gauss--Newton Joint Update}
\label{sec:damped_gn}

At iteration $k$, the local sensitivity matrix is estimated from the
performance and parameter changes between two consecutive iterates.
Specifically, for performance metric $p$ and optimized parameter $q$,
we use the successive finite-difference approximation
\begin{equation}
[\mathbf{J}_k]_{pq}
=
\frac{
P_{k,p}-P_{k-1,p}
}{
\phi_{k,q}-\phi_{k-1,q}
},
\qquad
p=1,\ldots,m,\quad q=1,\ldots,d .
\label{eq:secant_sensitivity}
\end{equation}

Because all parameters may change simultaneously between consecutive
iterates, Eq.~\eqref{eq:secant_sensitivity} represents a
trajectory-based numerical sensitivity rather than an independently
perturbed coordinate-wise partial derivative. The first sensitivity
estimate is initialized from two full-vector evaluations: the initial
configuration and a small simultaneous perturbation of all optimized
parameters. Thereafter, each newly evaluated HPO iterate is paired with
the preceding iterate to update the secant sensitivity, requiring no
additional coordinate-wise perturbation evaluations.

At iteration $k$, define the discrepancy between the desired and
observed performance vectors as
\begin{equation}
\mathbf{r}_k
=
\mathbf{P}^{\star}-\mathbf{P}_k .
\label{eq:target_residual}
\end{equation}
Using the resulting local sensitivity model, the desired joint parameter
change approximately satisfies
\begin{equation}
\mathbf{J}_k\Delta\boldsymbol{\phi}_k
\approx
\mathbf{r}_k .
\label{eq:desired_change}
\end{equation}

We compute the update using a Tikhonov-regularized damped
Gauss--Newton system,
\begin{equation}
\left(
\mathbf{J}_k^\top\mathbf{J}_k+
\lambda_d\mathbf{I}
\right)
\Delta\boldsymbol{\phi}_k
=
\mathbf{J}_k^\top\mathbf{r}_k ,
\label{eq:gn_system}
\end{equation}
or equivalently
\begin{equation}
\Delta\boldsymbol{\phi}_k
=
\left(
\mathbf{J}_k^\top\mathbf{J}_k+
\lambda_d\mathbf{I}
\right)^{-1}
\mathbf{J}_k^\top\mathbf{r}_k ,
\label{eq:gn_update}
\end{equation}
where $\lambda_d>0$ is the damping coefficient. The regularization
stabilizes the local inverse problem when $\mathbf{J}_k$ is
ill-conditioned or when $d>m$, and suppresses excessively large updates
caused by uncertain local sensitivity estimates.

A key property of Eq.~\eqref{eq:gn_system} is that it solves for
$\Delta\boldsymbol{\phi}_k\in\mathbb{R}^d$ as a single vector. The
sensitivity matrix therefore couples the measured performance responses
to all optimized parameters, and all $d$ coordinates are updated jointly
rather than through independent coordinate searches.

The next parameter vector is
\begin{equation}
\boldsymbol{\phi}_{k+1}
=
\Pi_{\Omega}
\left(
\boldsymbol{\phi}_k+
\alpha_k\Delta\boldsymbol{\phi}_k
\right),
\label{eq:parameter_update}
\end{equation}
where $\alpha_k$ is the update step length and $\Pi_{\Omega}$ restricts
the parameters to their prescribed feasible ranges when bounds are
imposed. In implementations using the adaptive step-length schedule,
\begin{equation}
\alpha_k
=
a+b/(1+k),
\label{eq:adaptive_step}
\end{equation}
where $a$ controls the asymptotic step length and $b$ controls the
additional early-iteration step size.

The complete updated vector $\boldsymbol{\phi}_{k+1}$ is then evaluated
once to obtain $\mathbf{P}_{k+1}$. Together with
$(\boldsymbol{\phi}_k,\mathbf{P}_k)$, this new observation provides the
successive finite-difference information used to construct
$\mathbf{J}_{k+1}$ through Eq.~\eqref{eq:secant_sensitivity}.

Consequently, following the two-evaluation full-vector sensitivity
initialization, each subsequent Gauss--Newton update introduces only one
new full-vector evaluation while updating all $d$ optimized parameters
simultaneously. This differs from a coordinate-wise finite-difference
Jacobian construction, which would require dedicated perturbation
evaluations for individual coordinates. When an evaluation entails
model training or post-training, the small local linear-algebra system
in Eq.~\eqref{eq:gn_system} is negligible compared with the model
evaluation cost.

\subsection{Training, Post-Training, and Post-Processing}
The same optimization loop is applied to training, post-training, and
post-processing, with candidate evaluations defined accordingly
(Table~\ref{tab:experiment_overview}). In all cases, successive
full-vector performance changes update the local sensitivity and all
hyperparameters are jointly updated while retaining the best observed
configuration.

\noindent
\begin{minipage}[t]{0.53\textwidth}
\vspace{0pt}

\begin{algorithm}[H]
\fontsize{7}{8}\selectfont
\caption{Iterative Gauss--Newton HPO}
\label{alg:newton_hpo}
\begin{algorithmic}[1]
\REQUIRE $\boldsymbol{\phi}_0$, target $\mathbf{P}^{\star}$,
$\lambda_d$, step rule $\alpha_k$, perturbation $\delta$, iterations $K$
\STATE Evaluate $\mathbf{P}_0=\mathbf{P}(\boldsymbol{\phi}_0)$
\STATE $\boldsymbol{\phi}_1\leftarrow
\boldsymbol{\phi}_0+\delta\mathbf{1}_d$;
evaluate $\mathbf{P}_1=\mathbf{P}(\boldsymbol{\phi}_1)$
\STATE Initialize secant sensitivity and best configuration
\FOR{$k=1,\ldots,K-1$}
    \STATE $\mathbf{r}_k \leftarrow \mathbf{P}^{\star}-\mathbf{P}_k$
    \STATE Update $\mathbf{J}_k$ from full-vector secant information
    \STATE Solve $(\mathbf{J}_k^\top\mathbf{J}_k+\lambda_d\mathbf{I})
    \Delta\boldsymbol{\phi}_k=\mathbf{J}_k^\top\mathbf{r}_k$
    \STATE $\boldsymbol{\phi}_{k+1}\leftarrow
    \Pi_{\Omega}(\boldsymbol{\phi}_k+\alpha_k\Delta\boldsymbol{\phi}_k)$
    \STATE Evaluate $\mathbf{P}_{k+1}=\mathbf{P}(\boldsymbol{\phi}_{k+1})$
    \STATE Update secant sensitivity and best configuration
\ENDFOR
\RETURN Best observed configuration
\end{algorithmic}
\end{algorithm}

\end{minipage}%
\hspace{0.02\textwidth}%
\begin{minipage}[t]{0.40\textwidth}
\vspace{0pt}
\centering

\includegraphics[width=0.88\linewidth]{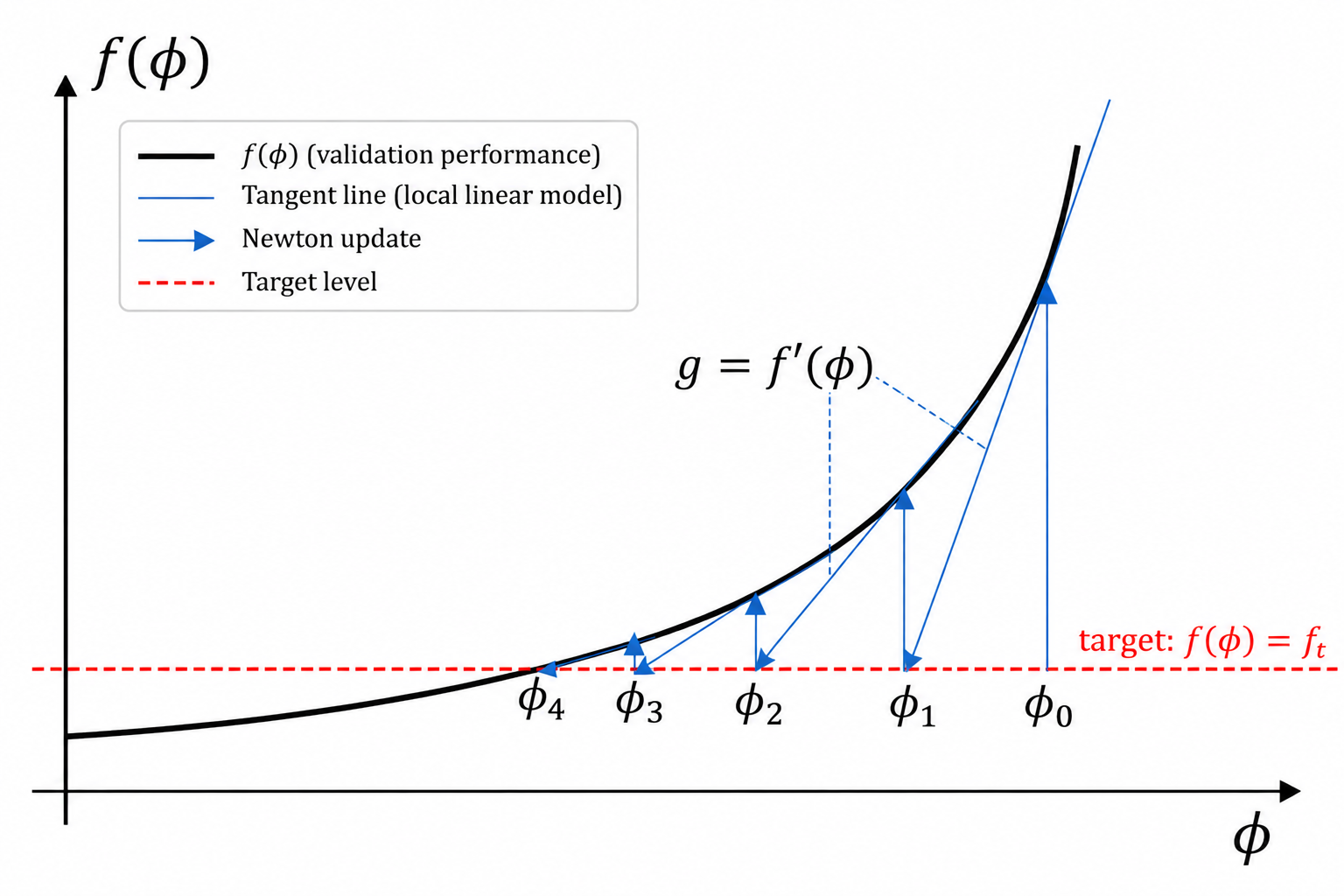}

\vspace{-1mm}
\captionof{figure}{\footnotesize
Scalar illustration of Newton-style target seeking; the
multi-parameter method jointly updates all hyperparameters.
}
\label{fig:newton}

\end{minipage}

\subsection{Experimental Design}
\label{sec:experimental_design}

We evaluate the framework in training-time HPO, VLM post-training, and
prediction-level post-processing. Table~\ref{tab:experiment_overview}
summarizes the corresponding parameter and response spaces and candidate
evaluations.

\begin{table*}[!ht]
\centering
\caption{Overview of the three experimental settings.}
\label{tab:experiment_overview}
\footnotesize
\setlength{\tabcolsep}{2pt}
\renewcommand{\arraystretch}{1.0}
\begin{tabular}{@{}lllll@{}}
\toprule
\textbf{Case} &
\textbf{ML setting} &
\textbf{Optimized parameters} &
\textbf{Candidate evaluation} &
\textbf{Responses} \\
\midrule
XGBoost &
Tabular classification &
4 training HPs &
Training + validation &
Accuracy \\
VLM &
Multimodal post-training &
6 SFT/LoRA HPs &
Restarted post-training + validation &
TPR, FPR \\
Threshold &
Classifier post-processing &
8 decision thresholds &
Fixed-score evaluation &
TP, FP \\
\bottomrule
\end{tabular}
\end{table*}

\section{Experiments}
\label{sec:experiments}
\subsection{Case I: Training-Time Hyperparameter Optimization with XGBoost}
\label{sec:xgboost}

\textbf{Experimental setup:}
We first evaluate the proposed method in a conventional tabular
classification setting, where each candidate evaluation requires
training a new model. The implementation code is publicly available in an anonymous
repository.\footnote{\url{https://github.com/vlm-grounding/newton_hp_search}} We use three public scikit-learn datasets:
Breast Cancer Wisconsin Diagnostic, Wine, and Digits, representing
binary, three-class, and ten-class classification, respectively.
For each candidate hyperparameter vector, a new XGBoost classifier is
trained and evaluated on a fixed validation split, and validation
accuracy provides the performance response used by the optimizer.
A held-out test split is excluded from hyperparameter selection.
Dataset characteristics, preprocessing, data splitting, parameter
ranges, fixed XGBoost settings, and baseline search configurations are
reported in Appendix~\ref{app:xgboost_details}.
We jointly optimize four continuous training hyperparameters, $\boldsymbol{\phi}
=[\eta,\;s,\;c,\;\lambda_{\mathrm{reg}}]^\top,$
corresponding to learning rate, row subsampling ratio, column
subsampling ratio, and $L_2$ regularization. All remaining XGBoost
settings are fixed across methods. Grid search, Random Search, and TPE
use the same hyperparameter domains and validation objective as the
proposed method. We report candidate evaluations explicitly so that
search budgets are compared in terms of expensive model
training--validation evaluations rather than optimizer iterations.

\textbf{Multi-dataset evaluation:}
We evaluate each method on Breast Cancer, Wine, and Digits across three
random seeds. Table~\ref{tab:xgb_results}(a) reports the mean and
standard deviation of the best observed validation accuracy. The
parenthesized values denote the number of candidate model evaluations
performed by each method.

Across the three datasets, the proposed method achieves validation
accuracy comparable to the black-box baselines under a smaller
candidate-evaluation budget, although it does not consistently attain
the highest accuracy (Table~\ref{tab:xgb_results}). These results show
that trajectory-based numerical sensitivity can produce competitive
configurations while jointly updating all four hyperparameters.

\begin{table*}[!ht]
\centering
\small

\caption{XGBoost evaluation. (a) Best validation accuracy across three
datasets (mean $\pm$ standard deviation over three seeds); parentheses
report the number of candidate model evaluations. (b) Sensitivity to
the Gauss--Newton damping coefficient on Breast Cancer (seed 42).}
\label{tab:xgb_results}

\vspace{2pt}

% ==================== LEFT ====================
\begin{minipage}[t]{0.62\textwidth}
\vspace{0pt}
\centering
\textbf{(a) XGBoost validation accuracy}\\[2pt]

\setlength{\tabcolsep}{3pt}
\begin{tabular}{lcccc}
\toprule
\textbf{Dataset} &
\shortstack{\textbf{Proposed}\\\textbf{(8)}} &
\shortstack{\textbf{Grid}\\\textbf{(320)}} &
\shortstack{\textbf{Random}\\\textbf{(20)}} &
\shortstack{\textbf{TPE}\\\textbf{(20)}} \\
\midrule
Breast Cancer & $.977{\pm}.011$ & $.976{\pm}.010$ & $.974{\pm}.007$ & $.977{\pm}.017$ \\
Wine          & $.991{\pm}.013$ & $.995{\pm}.006$ & $1.000{\pm}.000$ & $1.000{\pm}.000$ \\
Digits        & $.970{\pm}.011$ & $.976{\pm}.011$ & $.974{\pm}.012$ & $.977{\pm}.011$ \\
\bottomrule
\end{tabular}
\end{minipage}%
\hspace{0.01\textwidth}%
% ==================== RIGHT ====================
\begin{minipage}[t]{0.36\textwidth}
\vspace{0pt}
\centering
\textbf{(b) Effect of damping}\\[2pt]

\footnotesize
\setlength{\tabcolsep}{1.8pt}
\begin{tabular}{@{}lccccc@{}}
\toprule
$\lambda_d$ & $10^{-4}$ & $10^{-3}$ & $10^{-2}$ & $10^{-1}$ & $1$ \\
\midrule
Acc. & .9737 & .9737 & .9737 & .9649 & .9737 \\
\bottomrule
\end{tabular}

\end{minipage}
\end{table*}

\paragraph{Effect of damping:}
We additionally examine sensitivity to the Gauss--Newton damping
coefficient $\lambda_d$ on Breast Cancer using seed 42.
Table~\ref{tab:xgb_results}(b) shows relatively stable validation
accuracy across the tested range: four of the five damping values yield
$0.9737$, while $\lambda_d=10^{-1}$ yields $0.9649$.
Overall, the results indicate limited sensitivity to $\lambda_d$ over
the tested range.

\subsection{Case II: VLM Post-Training Hyperparameter Optimization}
\label{sec:vlm}

\textbf{Experimental setup:}
We next evaluate the proposed method in a substantially more expensive
multimodal post-training setting using Qwen3-VL-8B. In contrast to
Case I, where each hyperparameter vector controls a complete XGBoost
training run, here each candidate vector controls one LoRA/SFT
post-training run. Importantly, every candidate evaluation is restarted
from the same initial model and LoRA state, so that differences in
validation performance can be attributed to the post-training
hyperparameters rather than accumulated training across HPO iterations.
The resulting model is evaluated on the validation split, and the
observed TPR and FPR responses are used to construct the local
sensitivity estimate and determine the next joint hyperparameter update.
The held-out evaluation split is not used to construct the sensitivity
estimates or HPO updates and is reported only for evaluation.
The implementation code is publicly available in an anonymous
repository.\footnote{\url{https://github.com/vlm-grounding/Newton-VLM-PostTraining}}

We use a binary yes/no subset of VQAv2 with fixed training,
validation, and held-out splits of 5,000, 1,000, and 1,000 examples,
respectively. The held-out split is excluded from hyperparameter
optimization. Six continuous post-training hyperparameters are jointly
optimized: $\boldsymbol{\phi}=[\eta,\lambda_{\mathrm{wd}},p_{\mathrm{drop}},
\alpha_{\mathrm{LoRA}},\beta_1,\beta_2]^\top$, where $\eta$ is the
learning rate, $\lambda_{\mathrm{wd}}$ is weight decay,
$p_{\mathrm{drop}}$ is dropout, $\alpha_{\mathrm{LoRA}}$ is the LoRA
scaling parameter, and $\beta_1$ and $\beta_2$ are the AdamW momentum
coefficients.

TPR (recall) and FPR form the two-dimensional performance response used
by the target-seeking update. Precision is reported as a diagnostic
metric but is not used to construct the optimization update. Thus,
after sensitivity initialization, each HPO iteration jointly updates
all six hyperparameters, restarts LoRA/SFT post-training from the same
initial state, and performs one validation evaluation. Additional
implementation details and representative VQAv2 examples are provided
in Appendix~\ref{app:vlm_examples}.

\begin{figure}[!ht]
    \centering

    \begin{minipage}[t]{0.49\linewidth}
        \centering
        \includegraphics[
            width=\linewidth,
            height=3.5cm,
            keepaspectratio
        ]{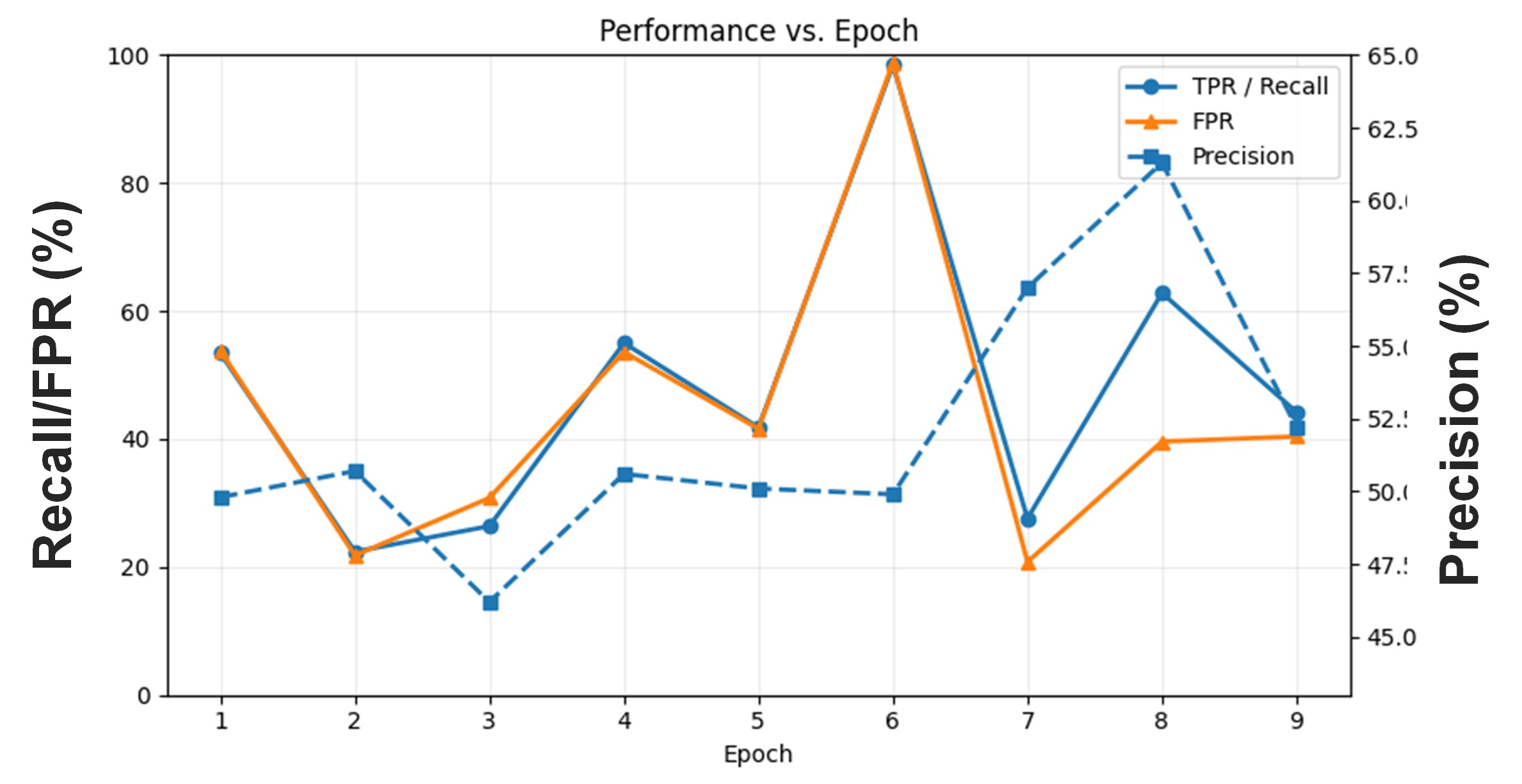}
        \vspace{-2mm}
        \centerline{\scriptsize (a) Validation performance}
    \end{minipage}
    \hfill
    \begin{minipage}[t]{0.49\linewidth}
        \centering
        \includegraphics[
            width=\linewidth,
            height=3.5cm,
            keepaspectratio
        ]{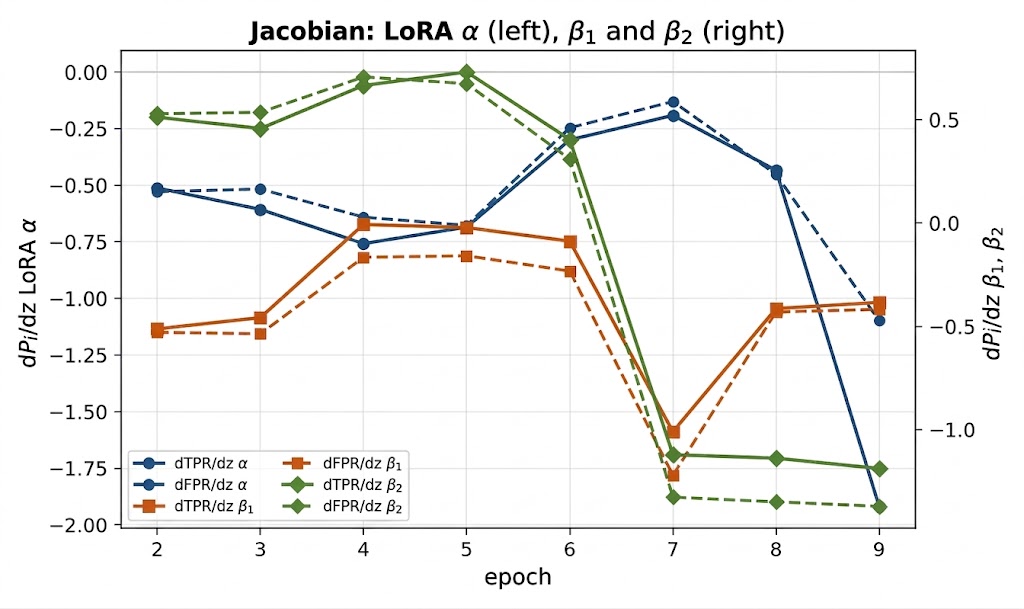}
        \vspace{-2mm}
        \centerline{\scriptsize (b) Local sensitivity}
    \end{minipage}

    \vspace{-1mm}
    \caption{VLM post-training across HPO iterations.
    (a) Validation TPR/FPR and diagnostic precision.
    (b) Local secant sensitivities of TPR/FPR with respect to
    $\alpha_{\mathrm{LoRA}}$, $\beta_1$, and $\beta_2$; operating-point
    changes coincide with changes in the estimated sensitivities.}
    \label{fig:vlm_performance}
    \vspace{-2mm}
\end{figure}

\textbf{Performance and local sensitivity evolution:}
Figure~\ref{fig:vlm_performance}(a) shows a strongly non-monotonic
validation trajectory across HPO iterations, while
Fig.~\ref{fig:vlm_performance}(b) shows the corresponding evolution
of the local secant sensitivities used by the optimization procedure.
The substantial variation in both the validation responses and their
estimated local sensitivities illustrates the trajectory-dependent
nature of the optimization landscape.

We additionally evaluate the checkpoints on the held-out split, which
is not used to construct the HPO updates. Held-out recall increases
from 53.1\% at iteration 1 to 54.6\% at iteration 4, while precision
increases from 48.7\% to 50.2\% and FPR decreases slightly from
54.2\% to 53.96\%.

These results also illustrate why recall alone is insufficient for this
target-seeking problem. Joint consideration of TPR and FPR distinguishes
increased recall accompanied by controlled false positives from increased
recall obtained primarily by shifting predictions toward the positive
class. Representative VQAv2 examples and prediction trajectories are
provided in Appendix~\ref{app:vlm_examples}. Individual predictions can
switch repeatedly between ``yes'' and ``no'', providing sample-level
examples of the non-monotonic behavior observed across HPO iterations.

\begin{figure*}[!ht]
    \centering
    \begin{minipage}[t]{0.49\textwidth}
        \centering
        \includegraphics[width=0.90\linewidth]
        {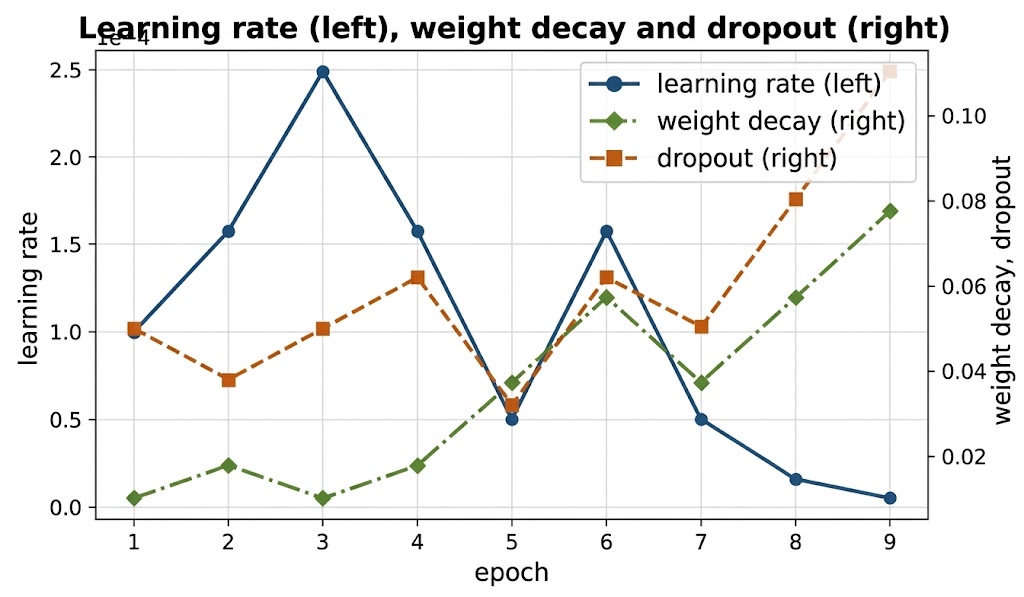}
        \vspace{-2mm}
        \centerline{\scriptsize (a) Learning rate, weight decay, and dropout}
    \end{minipage}
    \hfill
    \begin{minipage}[t]{0.49\textwidth}
        \centering
        \includegraphics[width=0.90\linewidth]
        {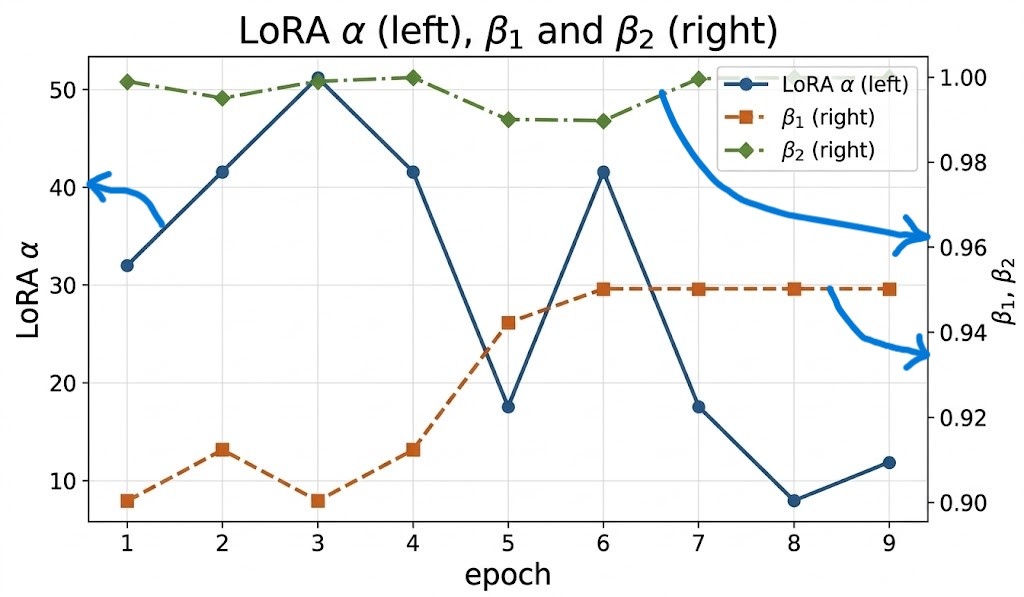}
        \vspace{-2mm}
        \centerline{\scriptsize (b) LoRA $\alpha$, $\beta_1$, and $\beta_2$}
    \end{minipage}

    \vspace{-1mm}
    \caption{Joint evolution of six VLM post-training hyperparameters.
    All parameters are updated simultaneously as one vector.}
    \label{fig:vlm_params}
    \vspace{-2mm}
\end{figure*}

\textbf{Joint hyperparameter evolution.}
Figure~\ref{fig:vlm_params} shows the corresponding trajectories of all
six post-training hyperparameters. Unlike coordinate-wise search, each
HPO step updates the complete six-dimensional hyperparameter vector
jointly.

The trajectories show nonuniform joint parameter changes. The learning
rate initially increases to approximately $2.5\times10^{-4}$ before
decreasing, while weight decay and dropout undergo smaller adjustments.
$\alpha_{\mathrm{LoRA}}$ varies more substantially in absolute scale,
whereas $\beta_1$ and $\beta_2$ remain within comparatively narrow
ranges. These different trajectories arise within the same full-vector
updates rather than from independent one-dimensional sweeps.

\subsection{Case III: Post-Processing Threshold Optimization}
\label{sec:threshold}

\textbf{Experimental design:}
We finally evaluate the proposed method in an eight-dimensional
post-processing problem requiring neither model training nor
post-training. A fixed proprietary CNN classifier produces probability
scores for eight anonymized classes. For class $i$, a prediction is
positive when $p_i\geq\phi_i$, with
$\boldsymbol{\phi}=[\phi_1,\ldots,\phi_8]^\top$ and
$0\leq\phi_i\leq1$. The evaluation set contains 6,487 instances,
including 1,002 positive and 5,485 negative cases. Model outputs remain
fixed throughout optimization; each candidate evaluation applies the
updated threshold vector to the stored scores and recomputes aggregate
TP and FP.

\begin{figure}[!ht]
    \centering
    \includegraphics[width=\linewidth]{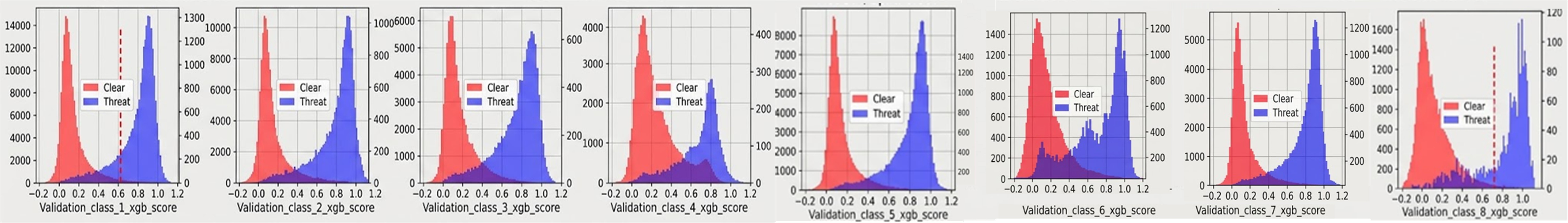}
    \caption{Prediction-score distributions for the eight-class
    threshold-optimization problem. Positive and negative scores exhibit
    different degrees of overlap across classes, motivating joint
    optimization of class-specific thresholds.}
    \label{fig:threshold_optimization}
\end{figure}

Figure~\ref{fig:threshold_optimization} shows representative positive
and negative score distributions. The degree of overlap varies across
classes, motivating class-specific thresholds, while the aggregate
TP/FP target requires the eight thresholds to be updated jointly.

\textbf{Target-seeking performance:}
Table~\ref{tab:newton_bayesian} compares the proposed method with a
Bayesian optimization reference for the 95/15 TP/FP target.

\begin{table}[!ht]
\centering
\caption{Comparison with Bayesian optimization for the 95/15 TP/FP
target. Target and outcome are reported in percent.}
\label{tab:newton_bayesian}
\scriptsize
\setlength{\tabcolsep}{2pt}
\resizebox{\columnwidth}{!}{%
\begin{tabular}{lcccc}
\toprule
\textbf{Method} &
\textbf{Input $\boldsymbol{\phi}$} &
\textbf{Target} &
\textbf{Outcome} &
\textbf{Output $\boldsymbol{\phi}$} \\
\midrule
Proposed &
$[0.5]^8$ &
95/15 &
94.7/16.3 &
$[.600,.599,.633,.701,.667,.598,.600,.701]$ \\
Bayesian$^{*}$ &
$\{[0.5]^8,[0.75]^8,[0.9]^8\}$ &
95/15 &
93.0/17.0 &
$[.567,.639,.433,.731,.575,.600,-1.23,2.57]$ \\
\bottomrule
\end{tabular}%
}
\vspace{1mm}
\parbox{\columnwidth}{\scriptsize
$^{*}$The Bayesian reference uses three initial design vectors and no
explicit $[0,1]$ box constraints; the infeasible outputs therefore
reflect the baseline configuration rather than an inherent limitation
of Bayesian optimization.}
\end{table}

The proposed method reaches TP/FP=94.7/16.3 at iteration 15, compared
with 93.0/17.0 for the Bayesian reference. Because the reference is
unconstrained and uses three initial design points, this result is
reported as a reference comparison rather than a controlled comparison
of the two algorithms.

\begin{figure*}[!ht]
\centering
\noindent
\begin{minipage}[t]{0.39\linewidth}
\vspace{0pt}
\centering

% ---- TABLE ----
\scriptsize
\setlength{\tabcolsep}{2pt}
\renewcommand{\arraystretch}{0.82}
\begin{tabular}{lccc}
\toprule
\textbf{Init. $\boldsymbol{\phi}_0$} &
\textbf{Target} &
\textbf{Outcome} &
\textbf{Iter.} \\
\midrule
$[.30]^8$ & 90/10 & 92.2/10.0 & $\sim$7 \\
$[.50]^8$ & 90/10 & 91.5/9.8  & 25 \\
$[.90]^8$ & 90/10 & 84.1/10.1 & 200 \\
$[.50]^8$ & 95/15 & 94.7/16.3 & 15 \\
Variable$^{*}$ & 96/15 & 94.0/15.3 & 15 \\
\bottomrule
\end{tabular}

\vspace{-1mm}
\captionof{table}{Effect of initialization and target operating point
(TP/FP, \%).}
\label{tab:threshold_seeds}

\end{minipage}%
\hspace{0.005\linewidth}%
\begin{minipage}[t]{0.595\linewidth}
\vspace{0pt}
\centering

% ---- FIGURE ----
\includegraphics[width=\linewidth]{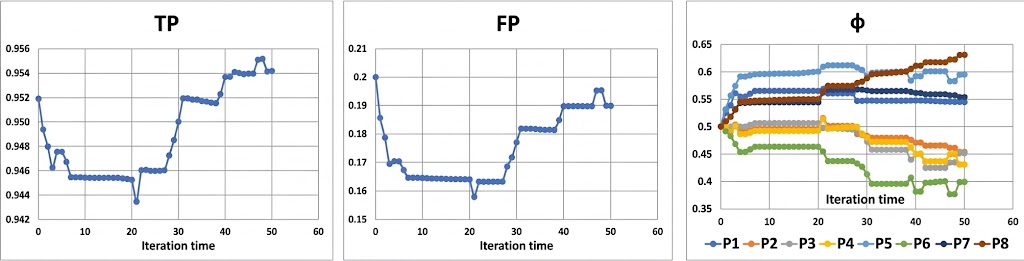}

\vspace{-1mm}
\captionof{figure}{
Trajectory for target TP$=95\%$, FP$=15\%$, initialized at
$[0.5]^8$. Best: 94.7/16.3 at iteration 15.
}
\label{fig:threshold_trajectory}

\end{minipage}
\end{figure*}

Table~\ref{tab:threshold_seeds} shows clear initialization dependence:
for the common 90/10 target, the three initializations follow different
trajectories and reach different operating points. This behavior is
consistent with the local, trajectory-dependent nature of the method
and motivates retaining the best observed iterate.

\textbf{Non-monotonic optimization dynamics:}
Figure~\ref{fig:threshold_trajectory} shows a representative trajectory
for the 95/15 target initialized from $[0.50]^8$. TP, FP, the threshold
parameters, and their locally estimated sensitivities vary
non-monotonically as the jointly updated thresholds move through
different local response regions. The non-monotonic trajectory shows that the best target-matching
operating point need not occur at the final iterate. The local
sensitivity is therefore repeatedly updated from successive
evaluations, while the best observed iterate is retained. Additional
trajectories for different initializations and targets are provided in
Appendix~\ref{app:threshold_demo}; they further illustrate the
initialization dependence and non-monotonic behavior of the local
search.
%%%%%%%%%%%%%%%%%%%%%%%%%%%%%%%%%%%%%%%%%%%%%%%%%%%%%%%%%%%%%%

\section{Conclusion}
\label{sec:conclusion}

We introduced a black-box, target-seeking HPO method that constructs
damped Gauss--Newton updates from successive full-vector evaluations,
without coordinate-wise perturbation runs or differentiation through
training. The method jointly updates all hyperparameters using local
secant sensitivities and Tikhonov regularization, while retaining the
best observed iterate to accommodate non-monotonic responses. Across
XGBoost training, Qwen3-VL-8B post-training, and classifier threshold
optimization, the experiments demonstrate the same update principle
under substantially different evaluation costs and response dimensions:
the method is competitive with standard HPO baselines on the tabular
benchmarks, supports joint six-dimensional VLM post-training updates
from TPR/FPR feedback, and approaches specified TP/FP operating points
in the eight-threshold problem. The observed non-monotonicity and
initialization dependence also highlight that the method is a local,
trajectory-dependent optimizer rather than a replacement for global
black-box HPO. Future work will investigate more robust sensitivity
estimation under noisy responses, larger model and search spaces, and
hybrid global--local optimization strategies.

\section{AI Usage}

AI-assisted tools were used to help identify relevant literature and
to support language editing and improvement of the manuscript. The
research methodology, experimental design, implementation, analysis,
and conclusions were developed and verified by the authors.

\bibliography{iclr2027_conference}
\bibliographystyle{iclr2027_conference}

\appendix

\section{Tabular Dataset and Experimental Details}
\label{app:xgboost_details}

Case~I evaluates XGBoost on three public datasets with different feature
representations and class structures. Breast Cancer Wisconsin Diagnostic
is a binary classification problem with 30 real-valued features describing
cell-nucleus morphology. Wine is a three-class classification problem
with 13 chemical measurements. Digits is a 10-class recognition problem
represented by 64 flattened $8\times8$ pixel-intensity features. Thus,
the experiments span binary and multiclass classification as well as
morphological, chemical, and image-derived feature representations.

\begin{table}[H]
\centering
\caption{Characteristics of the tabular classification tasks in Case~I.}
\label{tab:xgb_dataset_details}
\small
\setlength{\tabcolsep}{6pt}
\begin{tabular}{lrrl}
\toprule
\textbf{Dataset} & \textbf{Features} & \textbf{Classes} &
\textbf{Feature type} \\
\midrule
Breast Cancer & 30 & 2  & Nucleus morphology \\
Wine          & 13 & 3  & Chemical measurements \\
Digits        & 64 & 10 & $8\times8$ pixel intensities \\
\bottomrule
\end{tabular}
\end{table}

For each run, the data are stratified into training, validation, and
held-out test partitions. The proposed method and baseline HPO methods
use the same partition, with models fitted on the training set and
hyperparameters selected using validation performance. A median imputer
and standard scaler fitted on the training partition are subsequently
applied to the validation and test partitions without changing the
feature dimensionality. The test partition is excluded from the HPO
objective.

\section{VLM Examples and Prediction Variation}
\label{app:vlm_examples}

Case~II uses a selected binary yes/no subset of VQAv2 for VLM
post-training. The reported validation performance should therefore be
interpreted in the context of this selected subset rather than as the
overall VQAv2 performance of Qwen3-VL-8B. The subset contains
challenging examples requiring combinations of visual recognition,
semantic understanding, relational reasoning, and scene-level
interpretation rather than only direct object detection or simple
visual verification. The examples below provide qualitative context
for the relatively low absolute TPR and precision observed in
Section~\ref{sec:vlm}; the purpose of Case~II is to study HPO behavior
under a challenging multimodal post-training response rather than to
establish the benchmark-level capability of Qwen3-VL-8B.

\begin{figure}[H]
    \centering
    \includegraphics[width=\linewidth]{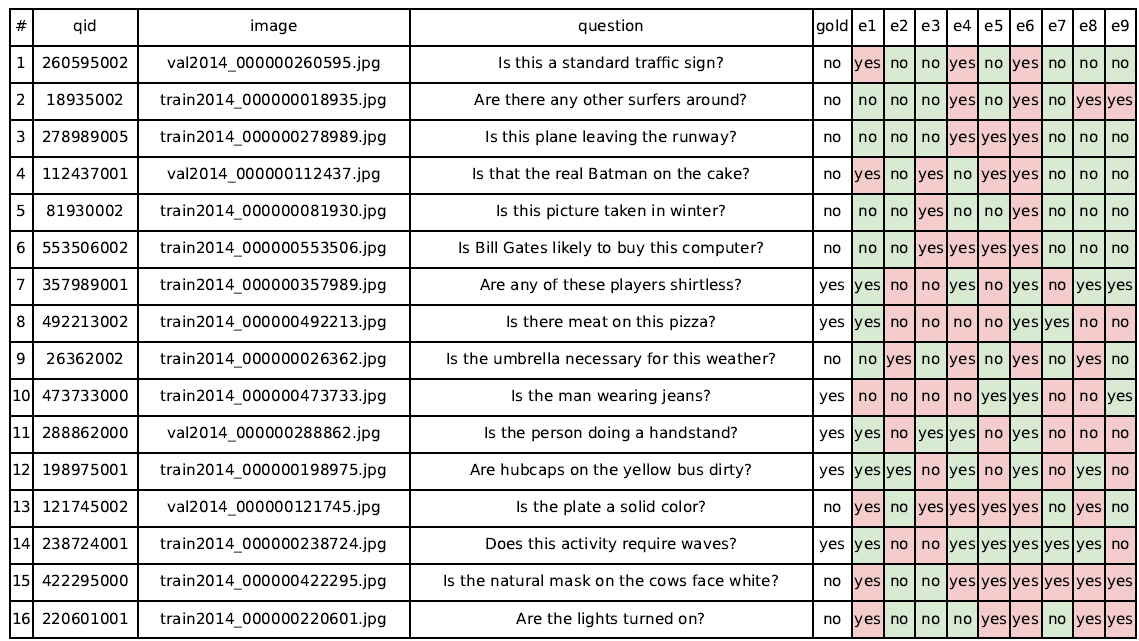}
    \caption{
    Predictions for 16 representative VQAv2 yes/no examples across
    successive VLM post-training epochs. The \texttt{gold} column
    denotes the ground-truth answer, and $e_1,\ldots,e_9$ show the
    corresponding model predictions at different epochs.
    }
    \label{fig:vqav2_epoch}
\end{figure}

Figure~\ref{fig:vqav2_epoch} presents 16 representative VQAv2 yes/no
questions and the corresponding model predictions across successive
post-training epochs. The \texttt{gold} column gives the ground-truth
answer, while $e_1,\ldots,e_9$ denote predictions obtained at different
epochs. Several examples switch repeatedly between ``yes'' and ``no'',
showing that sample-level prediction correctness does not necessarily
improve monotonically during post-training. This behavior provides
qualitative, sample-level context for the non-monotonic aggregate
validation responses observed in Case~II.

%~\ref{fig:vqa_image_examples}%
The following page shows the corresponding 16 input images together with their questions and ground-truth answers. The
examples span different visual and semantic requirements, including
object and attribute recognition, spatial and relational understanding,
scene interpretation, and questions requiring information to be
combined across multiple visual cues. They illustrate that the
evaluation subset contains more than simple object-presence questions
and provide additional context for the absolute performance observed
in the VLM experiment. Together with Figure~\ref{fig:vqav2_epoch},
these examples also show that prediction variation occurs across
different types of visual-question inputs rather than being confined
to a single example category.

% Include the qualitative VQA examples.
\includepdf[
    pages=-,
    pagecommand={\thispagestyle{plain}}
]{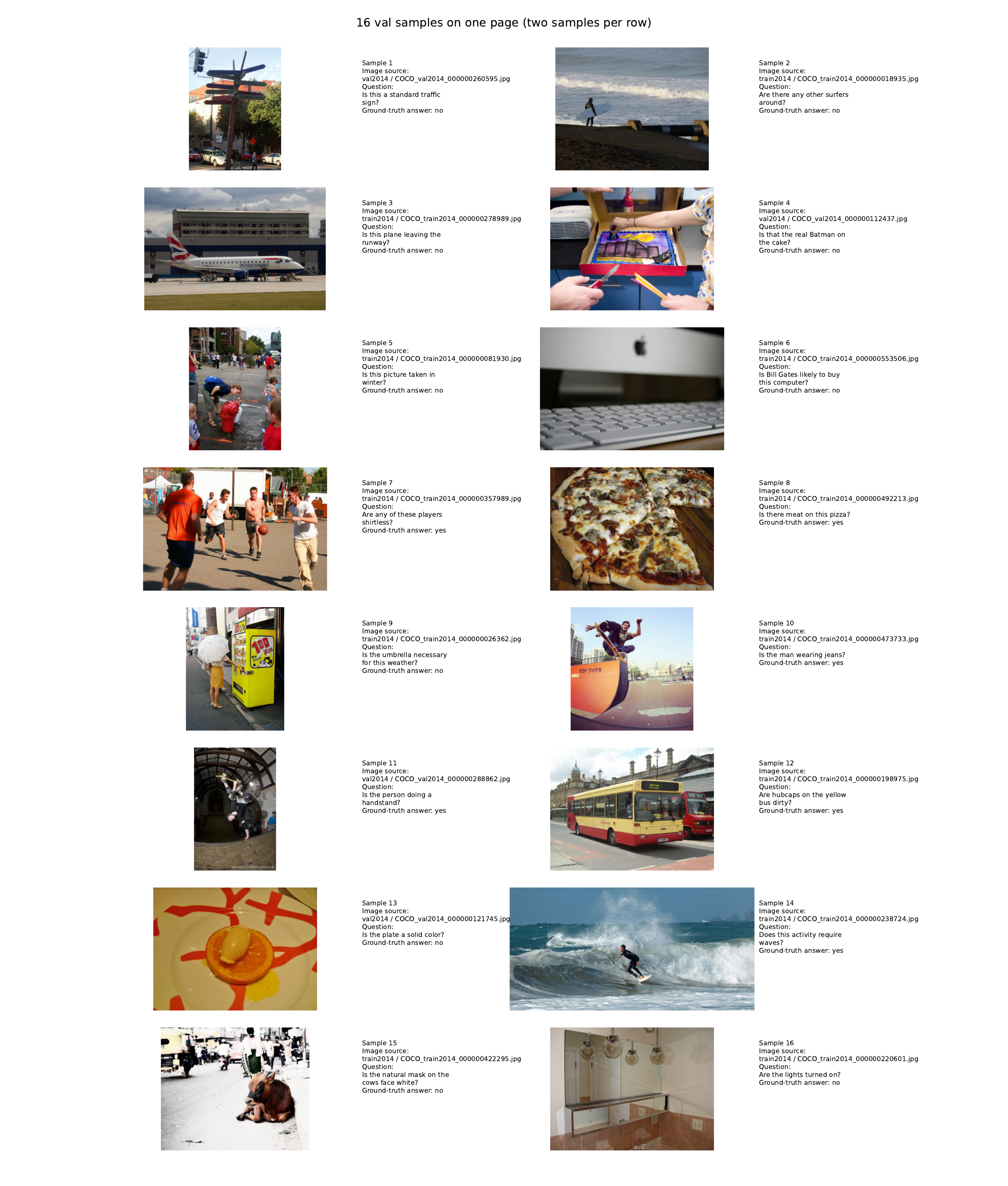}

\section{Additional Threshold-Optimization Trajectories}
\label{app:threshold_demo}

This section provides additional trajectories for the post-processing
threshold-optimization experiment in Section~\ref{sec:threshold}. All
experiments use the same fixed classifier outputs and jointly optimize
the eight class-specific probability thresholds. The initial threshold
vector and target TP/FP operating point are varied to examine their
effects on the local optimization trajectory. The representative 95/15
case is presented in Figure~\ref{fig:threshold_trajectory} of the main
text.

\begin{figure}[H]
\centering
\includegraphics[width=\linewidth]{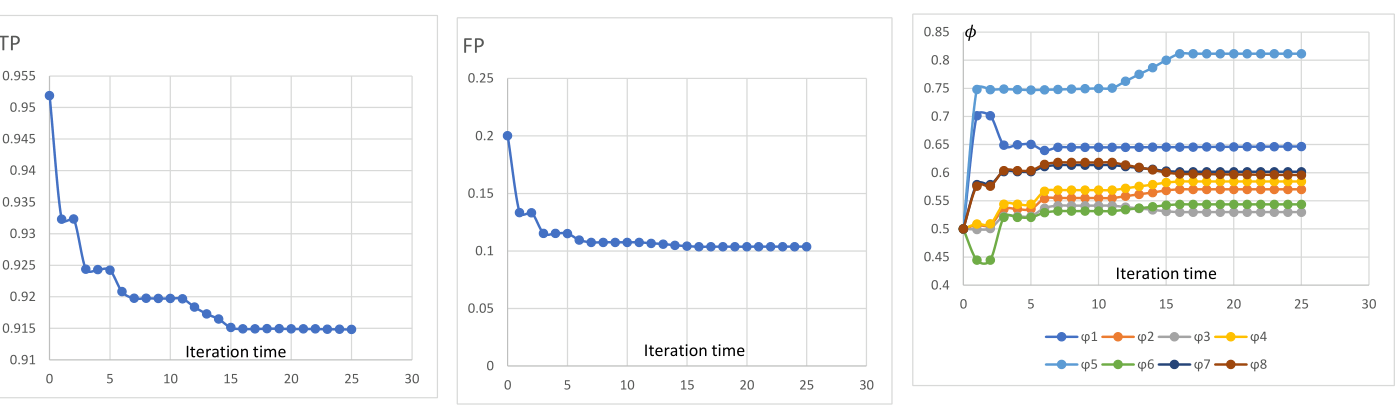}
\caption{
Additional trajectory for target TP$=90\%$ and FP$=10\%$, initialized
with $\boldsymbol{\phi}_0=[0.5]^8$. TP and FP decrease toward the target
region as the individual threshold coordinates evolve in different
directions. The best reported operating point is TP$=91.5\%$ and
FP$=9.8\%$ at iteration 25.
}
\label{fig:newton_seed05}
\end{figure}

\begin{figure}[H]
\centering
\includegraphics[width=\linewidth]{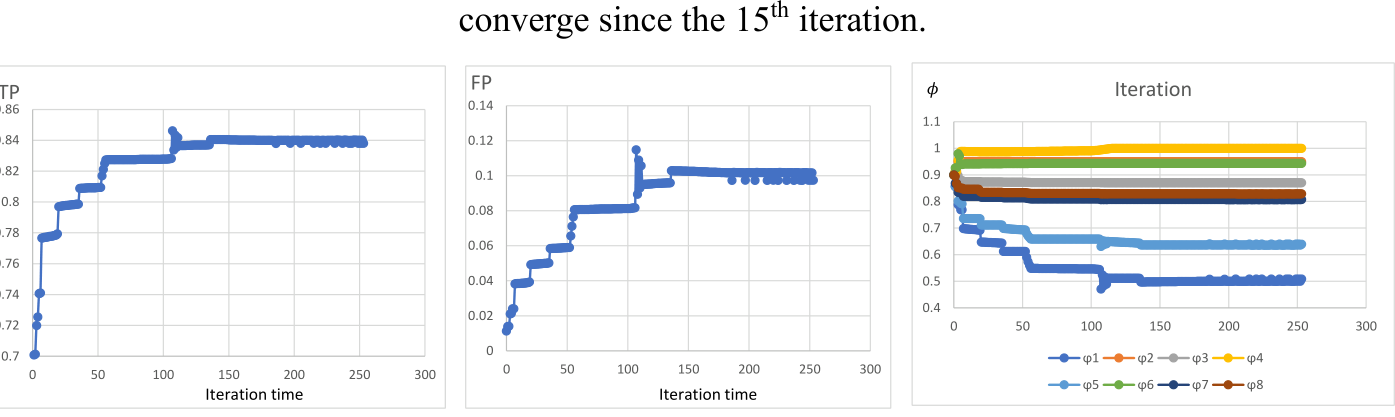}
\caption{
Additional trajectory for target TP$=90\%$ and FP$=10\%$, initialized
with the relatively high threshold vector
$\boldsymbol{\phi}_0=[0.9]^8$. TP and FP increase as the thresholds
move away from their high initial values. The reported operating point
at iteration 200 is TP$=84.1\%$ and FP$=10.1\%$, illustrating the
effect of initialization on the local optimization trajectory.
}
\label{fig:newton_seed09}
\end{figure}

\begin{figure}[H]
\centering
\includegraphics[width=\linewidth]{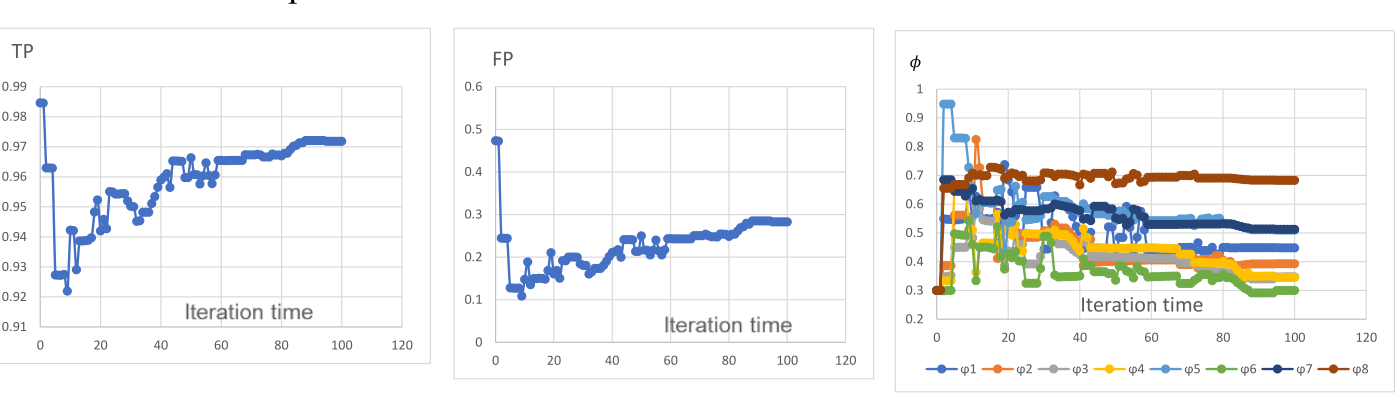}
\caption{
Additional trajectory for target TP$=90\%$ and FP$=10\%$, initialized
with $\boldsymbol{\phi}_0=[0.3]^8$. TP, FP, and individual threshold
values exhibit oscillatory behavior, with the best reported operating
point TP$=92.2\%$ and FP$=10.0\%$ occurring around iteration 7.
}
\label{fig:newton_seed03}
\end{figure}

\begin{figure}[t]
\centering
\includegraphics[width=\linewidth]{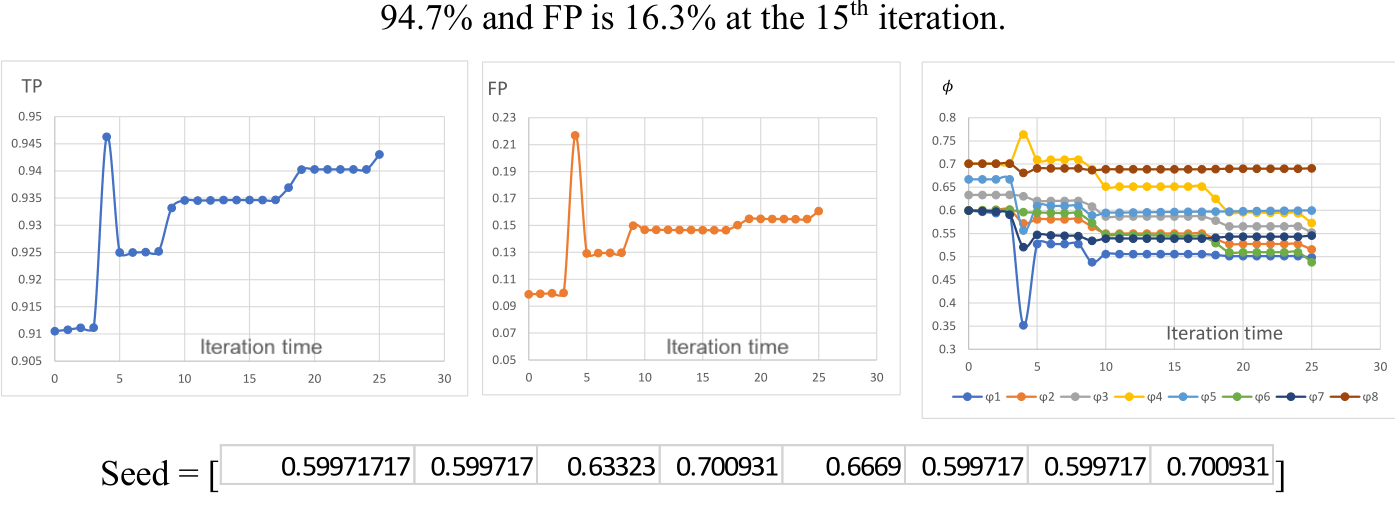}
\caption{
Additional trajectory for target TP$=96\%$ and FP$=15\%$ using the
nonuniform initialization
$\boldsymbol{\phi}_0=
[0.5997,0.5997,0.6332,0.7009,0.6669,0.5997,0.5997,0.7009]$,
obtained from preliminary optimization. TP and FP vary
non-monotonically during the trajectory, with the best reported
operating point TP$=94.0\%$ and FP$=15.3\%$ at iteration 15.
}
\label{fig:newton_variable_seed}
\end{figure}

Figures~\ref{fig:newton_seed05}--\ref{fig:newton_variable_seed} show
substantial differences in transient behavior across initializations.
For the common 90/10 target, initializations of $[0.3]^8$, $[0.5]^8$,
and $[0.9]^8$ produce different parameter paths, iteration horizons,
and achieved operating points. The trajectories also exhibit
non-monotonic or oscillatory changes in TP, FP, and individual
threshold coordinates. Together with the representative 95/15 case in
the main text, these results illustrate the initialization dependence
of the local search and motivate repeated re-estimation of the local
sensitivity and retention of the best observed iterate rather than
assuming monotonic convergence.

\end{document}